\documentclass{amsart}

\newtheorem{thm}{Theorem}[section]

\newtheorem{prop}[thm]{Proposition}

\newtheorem*{ggsconj}{Gerstenhaber-Giaquinto-Schack Conjecture}

\theoremstyle{definition}

\theoremstyle{remark}

\numberwithin{equation}{section}

\newcommand{\sln}{\mathfrak{sl}(n)}
\newcommand{\slf}{\mathfrak{sl}(5)}
\newcommand{\g}{\mathfrak{g}}
\newcommand{\frakh}{\mathfrak{h}}
\newcommand{\at}{\tilde{a}}
\newcommand{\qh}{\hat{q}}
\newcommand{\C}{\mathbb{C}}
\newcommand{\Z}{\mathbb{Z}}
\newcommand{\calt}{\mathcal{T}}

\DeclareMathOperator{\Tr}{Tr}

\begin{document}

%%%%%%%%%%%%%%%    TOPMATTER    %%%%%%%%%%%%%%%%

\title{Nonstandard solutions of the Yang-Baxter equation}

\author{Anthony Giaquinto}
\address{Department of Mathematics, Texas A \& M University,
College Station, TX 77843-3368}
\email{tonyg@math.tamu.edu}
\thanks{The first author was supported in part by a grant from the National Security 
Agency}
\author{Timothy J. Hodges} 
\address{University of Cincinnati, Cincinnati, OH 45221-0025,
U.S.A.}
\email{timothy.hodges@uc.edu}
\thanks{The second author was supported in part by a grant from the National 
Science Foundation}
\subjclass{Primary 81R50, 17B37; Secondary 16W30}

\date{November 17, 1997}

\keywords{Quantum group, R-matrix}

%%%%%%%%%%%%%%%   ABSTRACT    %%%%%%%%%%%%%

\begin{abstract}
    Explicit solutions of the quantum Yang-Baxter equation are given 
    corresponding to the 
    non-unitary solutions of the classical Yang-Baxter equation for $\slf$.
\end{abstract}

\maketitle

%%%%%%%%%%%%%%%% INTRODUCTION  %%%%%%%%%%%%

\section{Introduction}

	Etingof and Kazhdan recently proved that any finite dimensional 
	Lie bialgebra $\g$ may be quantized \cite{EK}. 
       That is, there exists a topological Hopf algebra structure on 
       $U(\g)[[h]]$ such that the Lie bialgebra 
       structure on $\g$ is the one induced on $\g$ by passing to the 
       ``semi-classical limit''. From this they 
       deduced a general procedure for quantizing solutions of the 
       classical Yang-Baxter equation (CYBE). 
       Thus, at least in theory, one can construct solutions of the quantum 
       Yang-Baxter equation from given 
       solutions of the classical Yang-Baxter equation. Unfortunately, 
       their procedure is not easy to 
       implement explicitly, even in small dimensional situations.

	In this note we exhibit an explicit answer to this problem for a 
	particularly interesting family 
	of Lie bialgebra structures on $\slf$. These are the bialgebra 
	structures associated to non-unitary 
	solutions of the CYBE (or equivalently of the modified 
	classical Yang-Baxter equation 
	(MCYBE)) as classified by Belavin and Drinfeld \cite{BD}. For each 
	such solution of the CYBE we 
	construct an $R$-matrix using the Gerstenhaber-Giaquinto-Schack (GGS) 
	conjecture \cite{GGS}. The YBE was 
	verified in each case using {\em Mathematica}.

	The GGS conjecture concerns the form of the quantization of such
	 solutions of the CYBE in the 
	case of $\sln$. The case of $\slf$ is to some extent the first 
	interesting case. For 
	$\mathfrak{sl}(2)$ there are no solutions of the MCYBE except the 
	standard one. For 
	$\mathfrak{sl}(3)$ the only non-standard solution is that associated
	to the
	well-known Cremmer-Gervais quantization and for $\mathfrak{sl}(4)$
	the nonstandard solutions are essentially of 
	three types, the 
	Cremmer-Gervais solution and 
       two other fairly simple 
	examples. The corresponding $R$-matrices for the latter two types can be 
	constructed using other techniques \cite{H4}. 
	On the other hand for $\slf$ there are 13 different types of 
	solutions to the MCYBE and for many of 
	these the corresponding $R$-matrix was hitherto unknown. The 
	validity of the GGS conjecture for 
	$\slf$ gives strong evidence that the conjecture should be 
	true for all $n$.

%%%%%%%%%%%%%%% MAIN BODY %%%%%%%%%%%%%%%%%
\section{Solutions to the CYBE and quantization}

 \subsection{The Belavin-Drinfeld description of solutions to the CYBE}
	Let $\g$ be a complex simple Lie algebra and let $\frakh$ be a 
	Cartan subalgebra. Let $\Delta$ be 
	the associated root system and $\Gamma$ a set of simple roots. 
	A classical $r$-matrix over $\g$ is 
	a an element $r \in \g \otimes \g$ satisfying the classical 
	Yang-Baxter equation
$$
[r_{12},r_{13}] + [r_{12},r_{23}]+[r_{13},r_{23}]=0
$$
Take an invariant bilinear from on form on $\g$ and let $t\in \g \otimes \g$ 
be the associated Casimir element. 
In \cite{BD} Belavin and Drinfeld gave the following description of solutions 
of the CYBE which satisfy 
$r_{12} +r_{21} =t$. These are the ``non-unitary'' solutions.

	Let $\Gamma_1$, $\Gamma_2$ be two subsets of $\Gamma$ and let 
	$\tau: \Gamma_1 \to \Gamma_2$ be a 
	bijection satisfying
\begin{enumerate}
\item $(\tau \alpha, \tau \beta)=(\alpha, \beta)$ for all
$\alpha,\beta 
\in \Gamma$;
\item For every $ \alpha \in \Gamma_1$, there is a $ k \geq 0$ with $\tau^k 
\alpha 
\in \Gamma_1 $ but 
$\tau^{k+1} \alpha \notin \Gamma_1 $.
\end{enumerate}
	The data $(\tau, \Gamma_1,\Gamma_2)$ (or more concisely just $\tau$) is often called a {\em Belavin-Drinfeld triple}. 
	Given such a triple $\tau$, an 
	element $r^0 \in \frakh \otimes \frakh$ is called $\tau$-admissible if 
\begin{enumerate}
\item $r^0_{12} +r^0_{21} =t^0$
\item $(\tau \alpha \otimes 1)r^0 + (1 \otimes \alpha)r^0 =t^0$
\end{enumerate}
where $t^0$ is the component of $t$ in $\frakh \otimes \frakh$.
A $\tau$-admissible $r^0$ is necessarily of the form 
$t^0/2 + \tilde{r}^0$ where $\tilde{r}^0\in \frakh \wedge \frakh$.
The set of all $\tilde{r}^0$ forms a linear subvariety of 
$\frakh \wedge \frakh $ of dimension $\binom {d}{2}$ where 
$d=\#(\Gamma - \Gamma_1)$.

Now $\tau$ can be extended to an isomorphism of Lie subalgebras 
$\tau : \g_1 \to \g_2$ where $\g_i$ is 
the Lie subalgebra of $\g$ associated to $\Gamma_i$. Choose 
$e_\alpha \in \g_\alpha$ such that 
$(e_\alpha, e_{-\alpha})=1$ and $\tau(e_\alpha) = e_{\tau \alpha}$ and define an
ordering on $\Delta$ by $\alpha \prec \beta$ if $\tau ^k\alpha =\beta$ for
some positive integer $k$. View $\g \wedge \g $ as a subset of $ \g \otimes \g$ 
via the identification $x\wedge y = 1/2(x\otimes y - y \otimes x)$.
Then Belavin and Drinfeld showed \cite{BD} 
that
$$
r = r^0 + \sum_{\alpha >0} e_{-\alpha} \otimes e_\alpha +
	 \sum_{\substack{\alpha,\beta > 0\\ \alpha \prec \beta}} e_{-\alpha} 
	\wedge e_\beta 
$$ 
is a solution of the Yang-Baxter equation satisfying $r_{12} +r_{21} =t$ 
and that every such solution is 
of this form for some choice of $\frakh$, $\Gamma$, $\tau$ and $r^0$.

For any $\g$ there is the ``trivial'' triple which has
$\Gamma_1 = \Gamma_2 =\emptyset$ and $\tilde{r}^0\in  \frakh \wedge \frakh$
arbitrary. A particularly interesting
triple for $\sln$ is the ``Cremmer-Gervais'' triple which has
$\Gamma_1=\{\alpha_2, \alpha_3, \ldots,\alpha_{n-1}\}, \quad 
	\Gamma_2=\{\alpha_1, \alpha_2,\ldots ,\alpha_{n-2}\}$, and 
	$\tau(\alpha_i) = \alpha_{i-1}$. 
In contrast to the trivial triple, there is a unique admissible $r^0$ for the
Cremmer-Gervais triple.

\subsection{The Gerstenhaber-Giaquinto-Schack conjecture}

	The Gerstenhaber-Gia\-quinto-Schack conjecture is a conjectured form 
	for the quantization of the 
	above classical $r$-matrices in the case where $\g = \mathfrak{sl}(n)$,
	considered as a subset of $M_{n}(\C)$. In this setting, a quantization
	of a classical $r$-matrix is an $R\in M_{n}(\C) \otimes M_{n}(\C)$
	which has semi-classical limit $r$ and satisfies the quantum
	Yang-Baxter equation $R_{12} R_{13}R_{23}=R_{23}R_{13}R_{12}.$
	
	Take the form to be the trace form $(x,y)=\Tr(xy)$ and let
	$\mathfrak h $ be the Cartan subalgebra consisting of diagonal
	matrices of trace zero. The standard Cartan-Weyl basis is then
	$e_{\alpha_i}=e_{i,i+1},\quad e_{-\alpha_i}=e_{i+1,i}$ and
	$h_{\alpha_i}=[e_{\alpha_i},e_{-\alpha_i}]=e_{ii}-e_{i+1,i+1}$.
	Let $\tau$ be a 
	Belavin-Drinfeld triple as described above and let 
	$r^0 \in \frakh \otimes \frakh$ be $\tau$-admissible. 
	Set 
$$a = \sum_{\substack{\alpha,\beta > 0\\ \alpha \prec \beta}} e_{-\alpha} 
\wedge e_\beta $$
and 
$$c_+ = \sum_{\alpha >0} e_{-\alpha} \otimes e_\alpha, \quad c = 
\sum_{\alpha >0} e_{-\alpha} \wedge e_\alpha.$$
 Set
$\epsilon = -(ac+ca+a^2)$. Now define $\at$ by
$$\at = \sum a^{ik}_{jl} q^{a^{ik}_{jl}\epsilon^{ik}_{jl}}e_{ij}\otimes e_{kl}$$
where $a = \sum a^{ik}_{jl} e_{ij}\otimes e_{kl}$ and similarly for $\epsilon$.
Set $\qh = q-q^{-1}$. The standard $R$-matrix is then  
$$R_s = q^{t^0+1/n} + \qh c_+ = q\sum _{i}e_{ii}\otimes e_{ii}+
\sum _{i\neq j}e_{ii}\otimes e_{jj}+ \qh \sum _{i>j}e_{ij}\otimes e_{ji}.$$ 
It is easy to check that $R_s$ satisfies the quantum Yang-Baxter
equation and that 
$PR_s$ satisfies the Hecke relation 
$(PR_s-q)(PR_s+q^{-1})=0$ where $P$ is the permutation matrix. 
\vspace{3mm}
\begin{ggsconj} Let $\tau $ be a Belavin-Drinfeld triple for $\sln$ and
suppose $r^0= t^0/2 + \tilde{r}^0$ is $\tau$-admissible. Then
the matrix
$$R= q^{\tilde{r}^0} (R_s +  \qh\, \at)\,q^{\tilde{r}^0} $$ 
satisfies the quantum Yang-Baxter equation and $PR$ satisfies
the Hecke relation.
\end{ggsconj}
\vspace{3mm}
Taking $\tau$ to be the
trivial triple yields the standard $R$-matrix 
 when $r^0=t^0/2$ and the standard multiparameter $R$-matrices when $r^0$ is 
 arbitrary. For use later, let $R(r^0)=q^{\tilde{r}^0} (R_s)q^{\tilde{r}^0}$ 
 denote the standard multiparameter $R$-matrix. 
 As is well known, if 
 $r^0 = t^0/2 + \sum_{i<j}c_{ij}\, e_{ii}\wedge e_{jj}$ then
 $$R(r^0)=q\sum _{i}e_{ii}\otimes e_{ii}+
\sum _{i<j}(q^{c_{ij}}e_{ii}\otimes e_{jj}+q^{-c_{ij}}e_{jj}\otimes e_{ii})+
\qh \sum _{i>j}e_{ij}\otimes e_{ji}.$$
 For the Cremmer-Gervais triples described above 
  the formula gives the 
  Cremmer-Gervais $R$-matrices 
  \cite{CG}. 
  
  \subsection{The GGS conjecture for $\slf$}
  We now consider the explicit form of the $R$-matrices associated to 
  the Belavin-Drinfeld triples on $\slf$. 
  According to the GGS Conjecture, each $R$ is of the form 
  $R(r^0) + \qh \, q^{\tilde{r}^0}\at q^{\tilde{r}^0}$ for an
  admissible $r^0$. The specific form of $R(r^0)$ has already
  been exhibited. The other summand, $q^{\tilde{r}^0}\at q^{\tilde{r}^0}$,
  is always a sum
  of ``quantized'' wedge products. Specifically, for positive roots
  $\alpha$ and $\beta$ and any constant $c$, 
  set $e_{-\alpha}\wedge_{c} e_{\beta}= q^{-c}e_{-\alpha}\otimes e_{\beta}-
  q^{c}e_{\beta}\otimes e_{-\alpha}$. For all triples, the term
  $q^{\tilde{r}^0}\at q^{\tilde{r}^0}$ is always of the form
  $\sum_{\substack{\alpha,\beta > 0\\ \alpha \prec \beta}}
  e_{-\alpha}\wedge_{c(\alpha,\beta)} e_{\beta} $ where the
  constants $c(\alpha,\beta)$ are determined by $\tilde{r}^0$ and $\epsilon$.
  
  Denote by $\calt$ the set of triples on $\slf$. Notice that if 
$(\tau, \Gamma_1,\Gamma_2)$ is a triple, 
then $(\tau^{-1}, \Gamma_2,\Gamma_1)$ is also a triple. Also the graph 
automorphism of $A_4$  induces a
bijection on the set of triples. Since these two involutions of $\calt$ commute, 
this gives an action of 
the group $\Z/2\Z\times \Z/2\Z$ on $\calt$.

\begin{prop} The Gerstenhaber-Giaquinto-Schack conjecture is true for $n=5$. The triples below comprise a complete set of
representatives 
from the $13$ orbits under 
the action of $C_2\times C_2$ on $\calt$. For each triple
the generic admissible $r^0$ and the Hecke $R$-matrix produced by the 
GGS conjecture are also explicitly given. 
\begin{enumerate}
\item $|\Gamma_1| = 3$
\vspace{5mm}
\begin{enumerate}
	\item The ``Cremmer-Gervais'' triple: 
	$\Gamma_1=\{\alpha_2, \alpha_3,\alpha_4\}$, \quad
	$\Gamma_2=\{\alpha_1, \alpha_2,\alpha_3\}$, 
	$\tau(\alpha_i) = \alpha_{i-1}$.
	\vspace{3mm}
     \begin{equation*}\begin{split}
	r^0&=t^0/2 + \frac{1}{5}( -3\,h_{\alpha_1}\wedge h_{\alpha_2} -
	4\,h_{\alpha_1}\wedge h_{\alpha_3} -3\, h_{\alpha_1}\wedge 
h_{\alpha_4}\\
	& \quad  -4\,h_{\alpha_2}\wedge h_{\alpha_3}
	 -4\,h_{\alpha_2}\wedge h_{\alpha_4}-3\,
	h_{\alpha_3}\wedge h_{\alpha_4}) \end{split}\end{equation*}
	\vspace{3mm}
	\begin{equation*}\begin{split}R&=R(r^0)+\qh \,(e_{54}\wedge_{2/5} e_{34} 
+
  e_{54}\wedge_{4/5} e_{23} 
    +e_{54}\wedge_{6/5}e_{12}
     + e_{43}\wedge_{2/5} e_{23}\\
    & \quad + e_{43}\wedge_{4/5} e_{12}+e_{32}\wedge_{2/5} e_{12}
    +e_{53}\wedge_{2/5} e_{24}
     +e_{53}\wedge_{4/5} e_{13}\\
   & \quad +e_{42}\wedge_{2/5} e_{13} +e_{52}\wedge_{2/5} e_{14})
   \end{split}\end{equation*}
	\vspace{5mm}
	\item The ``generalized Cremmer-Gervais" triple: 
	$\Gamma_1=\{\alpha_1, \alpha_3,\alpha_4\}$, \newline 
	$\Gamma_2=\{\alpha_1, \alpha_2,\alpha_4\}$, \quad 
	$\tau(\alpha_i) = \alpha_j$, where $j \equiv i+3$ (mod 5).
	\vspace{3mm}
     \begin{equation*}\begin{split}
	r^0&=t^0/2 + \frac{1}{5}( h_{\alpha_1}\wedge h_{\alpha_2} -
	2\,h_{\alpha_1}\wedge h_{\alpha_3} + h_{\alpha_1}\wedge h_{\alpha_4}\\
	 & \quad -2\,h_{\alpha_2}\wedge h_{\alpha_3} - 
	 2\,h_{\alpha_2}\wedge h_{\alpha_4}+
	h_{\alpha_3}\wedge h_{\alpha_4}) \end{split}\end{equation*}
	\vspace{3mm}
	\begin{equation*}\begin{split}R&=R(r^0)+\qh \,(e_{54}\wedge_{2/5} e_{23} 
+
  e_{21}\wedge_{4/5} e_{23} +
   e_{43}\wedge_{6/5}e_{23}
     + e_{21}\wedge_{2/5} e_{45}\\
     & \quad
    + e_{43}\wedge_{4/5} e_{45}+e_{43}\wedge_{2/5} e_{12}
    +e_{53}\wedge_{2/5} e_{13})\end{split}\end{equation*}
    \vspace{5mm}
    \end{enumerate}
\item $|\Gamma_1| = 2$
\vspace{5mm}
\begin{enumerate}
\item  $\Gamma_1=\{\alpha_3, \alpha_4\}$,\quad	
$\Gamma_2=\{\alpha_1, \alpha_2\}$,\quad	 $\tau(\alpha_i) = 
\alpha_{i-2}$
\vspace{3mm}
\begin{equation*}\begin{split}
r^0&=t^0/2 +  c\,h_{\alpha_1}\wedge h_{\alpha_2} +
       ((c-1)/2) \,h_{\alpha_1}\wedge h_{\alpha_3} 
	+c \,h_{\alpha_1}\wedge h_{\alpha_4} \\
	& \quad
	-((1+3c)/4)\,h_{\alpha_2}\wedge h_{\alpha_3} 
	+((c-1)/2)\,h_{\alpha_2}\wedge h_{\alpha_4} +
	c\,h_{\alpha_3}\wedge h_{\alpha_4} \end{split}\end{equation*}
	\vspace{3mm}
	$$R=R(r^0)+\qh\, (e_{43}\wedge_{(3+c)/8} e_{12} +
  e_{54}\wedge_{(3+c)/8} e_{23} +
   e_{53}\wedge_{(1-c)/2}e_{13})$$
   \vspace{5mm}
   \item  $\Gamma_1=\{\alpha_3, \alpha_4\}$,\quad		
$\Gamma_2=\{\alpha_1, \alpha_2\}$,\quad	 $\tau(\alpha_i) = 
\alpha_{5-i}$
\vspace{3mm}
\begin{equation*}\begin{split}
r^0&=t^0/2 + \frac{1}{5}( -c\,h_{\alpha_1}\wedge h_{\alpha_2}
       - (1+c)\, h_{\alpha_1}\wedge h_{\alpha_3} 
	-3 \,h_{\alpha_1}\wedge h_{\alpha_4} \\
	& \quad
	-2\,h_{\alpha_2}\wedge h_{\alpha_3} 
	+(c-1)\,h_{\alpha_2}\wedge h_{\alpha_4} +
	c\,h_{\alpha_3}\wedge h_{\alpha_4}) \end{split}\end{equation*}
   \vspace{3mm}
	$$R=R(r^0)+\qh \,(e_{54}\wedge_{3/5} e_{12} +
  e_{43}\wedge_{4/5} e_{23} +
   e_{53}\wedge_{-9/5}(-e_{13}))$$
   \vspace{5mm}
	\item	$\Gamma_1=\{\alpha_2, \alpha_4\}$, \quad	
$\Gamma_2=\{\alpha_1, \alpha_3\}$,\quad	  $\tau(\alpha_i) = \alpha_{i-1}$
\vspace{3mm}
\begin{equation*}\begin{split}
r^0&=t^0/2 +  c\,h_{\alpha_1}\wedge h_{\alpha_2}
       + (1+3c)\, h_{\alpha_1}\wedge h_{\alpha_3} 
	+(8c/3 + 1) \,h_{\alpha_1}\wedge h_{\alpha_4} \\
	& \quad
	(1+3c)\,h_{\alpha_2}\wedge h_{\alpha_3} 
	+(1+3c)\,h_{\alpha_2}\wedge h_{\alpha_4} +
	c \,h_{\alpha_3}\wedge h_{\alpha_4} \end{split}\end{equation*}
    \vspace{3mm}
	$$R=R(r^0)+\qh\, (e_{32}\wedge_{1+c} e_{12} +
   e_{54}\wedge_{1+c}e_{34})$$
   \vspace{5mm}
\item	$\Gamma_1=\{\alpha_2, \alpha_4\}$,\quad	 
$\Gamma_2=\{\alpha_1, \alpha_3\}$,\quad	  $\tau(\alpha_4) = \alpha_1$,\quad  
$\tau(\alpha_2)=\alpha_3$
\vspace{3mm}
\begin{equation*}\begin{split}
r^0&=t^0/2 + \frac{1}{5}( (2-c)\,h_{\alpha_1}\wedge h_{\alpha_2}
       + (1-c)\, h_{\alpha_1}\wedge h_{\alpha_3} 
	-3 \,h_{\alpha_1}\wedge h_{\alpha_4} \\
	& \quad
	+2\,h_{\alpha_2}\wedge h_{\alpha_3} 
	+(c-1)\,h_{\alpha_2}\wedge h_{\alpha_4} +
	c \,h_{\alpha_3}\wedge h_{\alpha_4}) \end{split}\end{equation*}
   \vspace{3mm}
$$R=R(r^0)+\qh \, (e_{54}\wedge_{2/5} e_{12} +
   e_{32}\wedge_{2/5}e_{34})$$
   \vspace{5mm}
\item $\Gamma_1=\{\alpha_1, \alpha_3\}$,\quad		
$\Gamma_2=\{\alpha_1, \alpha_4\}$,\quad	  $\tau(\alpha_i) = \alpha_j$,
 where $j \equiv i+3$ (mod 5).
 \vspace{3mm}
 \begin{equation*}\begin{split}
r^0&=t^0/2 +  ((1-3c)/2)\,h_{\alpha_1}\wedge h_{\alpha_2}
       + ((c-1)/2)\, h_{\alpha_1}\wedge h_{\alpha_3} 
	+c \,h_{\alpha_1}\wedge h_{\alpha_4} \\
	& \quad
	+(3c-1)\,h_{\alpha_2}\wedge h_{\alpha_3} 
	+(3c-1)\,h_{\alpha_2}\wedge h_{\alpha_4} +
	c \,h_{\alpha_3}\wedge h_{\alpha_4} \end{split}\end{equation*}
	\vspace{3mm}
	$$R=R(r^0)+\qh\, (e_{43}\wedge_{(1+3c)/4} e_{12} +
   e_{21}\wedge_{(1+3c)/4}e_{45}+
   e_{43}\wedge_{1-c}e_{45})$$
   \vspace{5mm}
 \item $\Gamma_1=\{\alpha_3, \alpha_4\}$,\quad		
$\Gamma_2=\{\alpha_1, \alpha_2\}$, \quad  $\tau(\alpha_i) = 
\alpha_{i-1}$
\vspace{3mm}
\begin{equation*}\begin{split}
r^0&=t^0/2 +  ((c-3)/6)\,h_{\alpha_1}\wedge h_{\alpha_2}
       + ((c-1)/2)\, h_{\alpha_1}\wedge h_{\alpha_3} 
	+c\, h_{\alpha_1}\wedge h_{\alpha_4} \\
	& \quad
	+((c-1)/2)\,h_{\alpha_2}\wedge h_{\alpha_3} 
	+(4c/3)\,h_{\alpha_2}\wedge h_{\alpha_4} +
	c \,h_{\alpha_3}\wedge h_{\alpha_4} \end{split}\end{equation*}
	\vspace{3mm}
	$$R=R(r^0)+\qh \, (e_{32}\wedge_{1/2+c/6} e_{12} +
   e_{43}\wedge_{1/2+c/6}e_{23}+
   e_{43}\wedge_{1+c/3}e_{12})$$
   \vspace{5mm}
 \end{enumerate}

\item $|\Gamma_1| = 1$
\vspace{5mm}
\begin{enumerate}
\item	 $\Gamma_1=\{\alpha_1 \}$,\quad	 $\Gamma_2=\{\alpha_2 \}$, \quad	
$\tau(\alpha_1) = \alpha_{2}$
\vspace{3mm}
\begin{equation*}\begin{split}
r^0&=t^0/2 + ((1+y)/3)\, h_{\alpha_1}\wedge h_{\alpha_2}
       + y\, h_{\alpha_1}\wedge h_{\alpha_3} 
	+((3z-x)/3)\, h_{\alpha_1}\wedge h_{\alpha_4}\\ 
	& \quad +y\,h_{\alpha_2}\wedge h_{\alpha_3} 
	+z\,h_{\alpha_2}\wedge h_{\alpha_4} +
	x\, h_{\alpha_3}\wedge h_{\alpha_4} 
	\end{split}\end{equation*}
	\vspace{3mm}
	$$R=R(r^0)+\qh \, (e_{21}\wedge_{(2-y)/3} e_{23} )$$
	\vspace{5mm}
	\item	 $\Gamma_1=\{\alpha_1 \}$,\quad	 $\Gamma_2=\{\alpha_3 \}$,\quad	 
$\tau(\alpha_1) = \alpha_{3}$
\vspace{3mm}
\begin{equation*}\begin{split}
r^0&=t^0/2 + ((z-2y)/2)\, h_{\alpha_1}\wedge h_{\alpha_2}
       + ((1+x)/2)\, h_{\alpha_1}\wedge h_{\alpha_3} 
	 +x\, h_{\alpha_1}\wedge h_{\alpha_4}\\ 
	& \quad +y\,h_{\alpha_2}\wedge h_{\alpha_3}
	+z\,h_{\alpha_2}\wedge h_{\alpha_4} +
	x\, h_{\alpha_3}\wedge h_{\alpha_4} \end{split}\end{equation*}
	\vspace{3mm}
	$$R=R(r^0)+\qh\, (e_{21}\wedge_{(z-2)/4} e_{34} )$$
	\vspace{5mm}
	\item	 $\Gamma_1=\{\alpha_1 \}$,\quad	 $\Gamma_2=\{\alpha_4 \}$,\quad	 
$\tau(\alpha_1) = \alpha_{4}$
\begin{equation*}\begin{split}
\vspace{3mm}
r^0&=t^0/2 + ((y-2z)/2)\,h_{\alpha_1}\wedge h_{\alpha_2}
       + ((y-2x)/2) h_{\alpha_1}\wedge h_{\alpha_3} 
	+((1-x+z)/2)\, h_{\alpha_1}\wedge h_{\alpha_4}\\ 
	& \quad + y\, h_{\alpha_2}\wedge h_{\alpha_3}
	+ z\,h_{\alpha_2}\wedge h_{\alpha_4} +
	x \,h_{\alpha_3}\wedge h_{\alpha_4} \end{split}\end{equation*}
	\vspace{3mm}
	$$R=R(r^0)+\qh \,(e_{21}\wedge_{(y+2)/4} e_{45} )$$
	\vspace{5mm}
	\item	 $\Gamma_1=\{\alpha_2 \}$,\quad	 $\Gamma_2=\{\alpha_3 \}$,\quad	 
$\tau(\alpha_2) = \alpha_{3}$
 \begin{equation*}\begin{split}
 \vspace{3mm}
r^0&=t^0/2 + (-1+3y-z)\,h_{\alpha_1}\wedge h_{\alpha_2}
       + (-1-x+3y) h_{\alpha_1}\wedge h_{\alpha_3} +
	3(z-x)\, h_{\alpha_1}\wedge h_{\alpha_4}\\ 
	& \quad +y\, h_{\alpha_2}\wedge h_{\alpha_3} 
	+z\,h_{\alpha_2}\wedge h_{\alpha_4} +
	x\, h_{\alpha_3}\wedge h_{\alpha_4} \end{split}\end{equation*}
	\vspace{3mm}
	$$R=R(r^0)+\qh \,(e_{32}\wedge_{1-x-y+z} e_{34} )$$
	\end{enumerate}
	\vspace{5mm}
	\item $|\Gamma_1| = 0$
	\vspace{5mm}
	The ``trivial triple:'' $\Gamma_1 =\Gamma_2 =\emptyset$\\
	\vspace{3mm}
	$r^0 = t^0/2+ \tilde{r}^0$ with $\tilde{r}^0\in \mathfrak h \wedge
	\mathfrak h$ arbitrary.\\
	 \vspace{3mm}
	 $R = R(r^0)$ is the standard multiparameter $R$-matrix.
  \end{enumerate}
  \end{prop}
Perhaps the most interesting new $R$-matrix is that associated to type
$1\,(b)$, the generalized Cremmer-Gervais triple. Like the
Cremmer-Gervais triple, its $\Gamma_1$, which must omit at least one root,
omits precisely one and thus its $r^0$ is uniquely determined. 
Setting $\hat{p}=-\qh$, the matrix form of the generalized Cremmer-Gervais
$R$-matrix is
    \vspace{2mm}
    $$\left(
\begin{smallmatrix} 
q& 0& 0& 0& 0& 0& 0& 0& 0& 0& 0& 0& 0& 0& 0& 0& 0& 0& 0& 0& 0& 0& 0& 0& 0\\
0& q^{\frac{1}{5}}& 0& 0& 0& 0& 0& 0& 0& 0& 0& 0& 0& 0& 0& 0& 0& 
        0& 0& 0& 0& 0& 0& 0& 0\\
   0& 0& q^{-\frac{3}{5}}& 0& 0& 0& 0& 0& 0& 0& 0& 0& 0& 0& 0& 0& 0& 0& 0& 0& 
   0& 0& 0& 0& 0\\
        0& 0& 0& q^{\frac{3}{5}}& 
        0& 0& 0& \hat{p}  q^{\frac{2}{5}}& 0& 0& 0&
        0& 0& 0& 0& 0& 0& 0& 0& 
        0& 0& 0& 0& 0& 0\\
        0& 0& 0& 0& q^{-\frac{1}{5}}& 0& 0& 0& 0& 0& 0& 0& 
        \hat{p}  q^{\frac{2}{5}}& 0& 0& 0& 0& 0& 0& 0& 0& 0& 0& 0& 0\\       
        0& 
        \hat{q}& 0& 0& 0& q^{-\frac{1}{5}}& 0& 0& 0& 0& 0& 0& 0& 0& 0& 
        0& 0& 0& 0& 0& 0& 0& 0& 0& 0\\
         0& 0& 
        \hat{q}q^{-\frac{4}{5}}  &      
         0& 0& 0& q& 0& 0& 0& 
        \hat{p} q^{\frac{4}{5}}& 0& 0& 0& 0& 0& 0& 0& 0& 0& 0& 0& 0& 0& 
        0\\       
        0& 0& 0& 0& 0& 0& 0& q^{\frac{1}{5}}& 0& 0& 0& 0& 0& 0& 0& 0& 0& 0& 
        0& 0& 0& 0& 0& 0& 0\\        
        0& 0& 0& 0& \hat{q}q^{-\frac{2}{5}}& 
        0& 0& 0& q^{-\frac{3}{5}}& 0& 0& 0& \hat{p}q^{\frac{6}{5}}& 0& 0& 0& 
        0& 0& 0& 0& 0& 0& 0& 0& 0\\        
        0& 0& 0& 0& 0& 0& 0& 0& 0& q^{\frac{3}{5}}& 
        0& 0& 0& \hat{p}q^{\frac{2}{5}}& 0& 0& 0& 0& 0& 0& 0& 0& 0& 0& 
        0\\        
        0& 0& \hat{q}& 0& 0& 0& 0& 0& 0& 0& q^{\frac{3}{5}}& 0& 0& 
        0& 0& 0& 0& 0& 0& 0& 0& 0& 0& 0& 0\\
        0& 0& 0& 0& 0& 0& 0& 
        \hat{q}& 0& 0& 0& q^{-\frac{1}{5}}& 0& 0& 0& 0& 0& 0& 0& 0& 0& 
        0& 0& 0& 0\\       
        0& 0& 0& 0& 0& 0& 0& 0& 0& 0& 0& 0& q& 0& 0& 0& 0& 0& 
        0& 0& 0& 0& 0& 0& 0\\     
        0& 0& 0& 0& 0& 0& 0& 0& 0& 0& 0& 0& 0& 
        q^{\frac{1}{5}}& 0& 0& 0& 0& 0& 0& 0& 0& 0& 0& 0\\       
       0& 0& 0& 0& 0& 0& 0& 
        0& 0& 0& 0& 0& 0& 0& q^{-\frac{3}{5}}& 0& 0& 0& 0& 0& 0& 0& 0& 0& 0\\
        0& 0& 0& \hat{q}& 0& 0& 0& 0& 0& 0& 0& 
        \hat{q}q^{-\frac{2}{5}}& 0& 0& 0& q^{-\frac{3}{5}}& 0& 0& 0& 
        0& 0& 0& 0& 0& 0\\
        0& 0& 0& 0& 0& 0& 0& 0& \hat{q}& 0& 0& 
        0& \hat{q}q^{-\frac{6}{5}}& 0& 0& 0& q^{\frac{3}{5}}& 0& 0& 0& 
        \hat{p} q^{\frac{2}{5}}& 0& 0& 0& 0\\
         0& 0& 0& 0& 0& 0& 0& 0& 0& 
        0& 0& 0& 0& \hat{q}& 0& 0& 0& q^{-\frac{1}{5}}& 0& 0& 0& 0& 0& 
        0& 0\\
         0& 0& 0& 0& 0& 0& 0& 0& 0& 0& 0& 0& 0& 0& 
        \hat{q}q^{-\frac{4}{5}}& 0& 0& 0& q& 0& 0& 0& 
         \hat{p}q^{\frac{4}{5}}& 0& 0\\
         0& 0& 0& 0& 0& 0& 0& 0& 0& 0& 0& 
        0& 0& 0& 0& 0& 0& 0& 0& q^{\frac{1}{5}}& 0& 0& 0& 0& 0\\
         0& 0& 0& 0& 
        \hat{q}& 0& 0& 0& 0& 0& 0& 0& 
        \hat{q}q^{-\frac{2}{5}}& 0& 0& 0& 0& 0& 0& 0& q^{\frac{1}{5}}& 
        0& 0& 0& 0\\
        0& 0& 0& 0& 0& 0& 0& 0& 0& \hat{q}& 0& 0& 0& 
        0& 0& 0& 0& \hat{q}q^{-\frac{2}{5}}& 0& 0& 0& 
        q^{-\frac{3}{5}}& 0& 0& 0\\
        0& 0& 0& 0& 0& 0& 0& 0& 0& 0& 0& 0& 0& 0& 
        \hat{q}& 0& 0& 0& 0& 0& 0& 0& q^{\frac{3}{5}}& 0& 0\\
        0& 0& 0& 
        0& 0& 0& 0& 0& 0& 0& 0& 0& 0& 0& 0& 0& 0& 0& 0& \hat{q}& 0& 
        0& 0& q^{-\frac{1}{5}}& 0\\
        0& 0& 0& 0& 0& 0& 0& 0& 0& 0& 0& 0& 0& 0& 
        0& 0& 0& 0& 0& 0& 0& 0& 0& 0& q \end{smallmatrix}\right). $$
	
	\vspace{5mm}

%%%%%%%%%  CONCLUSION %%%%%%%%%%%%%%%
\section{Conclusion}

	We have constructed here quantizations of each type of 
	non-unitary solution of the classical Yang-Baxter equation for 
	$\slf$. In so doing we verified in this case the conjecture of 
	Gerstenhaber, Giaquinto and Schack. 
	This gives further evidence that the GGS conjecture should be true
	for all Belavin-Drinfeld triples on $\sln$.

	One can proceed in the usual way to construct for each of these $R$, a quantization of 
	$\C[SL(5)]$, the algebra of algebraic functions on $SL_\C(5)$. First one constructs the associated 
	bialgebra $A(R)$. Using a case-by-case analysis one can see that 
	the Poincare series of the 
	associated quantum space and exterior algebra are the same as in the 
	commutative case. Thus $A(R)$ contains a group-like $q$-determinant 
	element $D$ which turns out to be central. Hence one may define a 
	Hopf algebra structure on $\C_R[SL(5)] =A(R)/(D-1)$. 
	Since $R$ is a Hecke symmetry in the sense
	of Gurevich, it is possible to exploit some
	Hecke algebra techniques to show that the category of 
	comodules over these Hopf algebras is equivalent as a rigid monoidal 
	category to the category of comodules 
	over $\C_q[SL(5)]$ \cite{Gu,Hai,KW}. 
	Hence these $R$-matrices do produce genuine nonstandard quantizations 
	of $\C[SL(5)]$.

%%%%%%%%   REFERENCES %%%%%%%%%%%%%%%

\end{document}